\newcommand{\CLASSINPUTinnersidemargin}{12mm}
\newcommand{\CLASSINPUToutersidemargin}{12mm}
\documentclass[10pt, final, journal, letterpaper, twoside, column]{IEEEtran}
\usepackage{graphicx}
\usepackage{subcaption}
\usepackage{epsfig} 
\usepackage{mathptmx} 
\usepackage{times} 
\usepackage{amsmath} 
\usepackage{amssymb}
\usepackage{enumerate}
\usepackage{mathrsfs}
\usepackage{euscript}
\usepackage{psfrag}
\usepackage[section]{placeins} 
\usepackage{amsbsy}
\usepackage{textcomp}
\usepackage{blkarray}
\usepackage{mathdots}
\usepackage[usenames,dvipsnames]{color}
\usepackage{mdframed}
\usepackage{caption}
\usepackage{subcaption}
\usepackage[T1]{fontenc}

\newtheorem{theorem}{Theorem}[section]

\newtheorem{corollary}[theorem]{Corollary}
\newtheorem{lemma}[theorem]{Lemma}
\newtheorem{example}[theorem]{Example}
\newtheorem{remark}[theorem]{Remark}

\begin{document}
  \title{Geometry and control of the nonholonomic integrator: An electrodynamics analogy}
 \author{Pragada Shivaramakrishna, A. Sanand Amita Dilip \thanks{Pragada Shivaramakrishna is with the Department of Aerospace Engineering, 
 Indian Institute of Technology Kharagpur, India. {\tt\small shivaramkratos@gmail.com, shivaram@iitkgp.ac.in}}\thanks{
  A. Sanand Amita Dilip is with the Department of Electrical Engineering, Indian Institute of Technology Kharagpur, India. 
 {\tt\small sanand@ee.iitkgp.ac.in}} 
}
 \maketitle
 \begin{abstract}
  We consider some generalizations of the classical nonholonomic integrator and give a geometric approach to characterize controllability for these systems. We use Stokes' theorem 
and results from complex analysis to obtain necessary and sufficient conditions for controllability of these systems. Furthermore, we show that optimal 
trajectories of certain minimum energy optimal control problems defined on these systems can be identified with the trajectory of a charged particle in an 
electromagnetic field. 

 \end{abstract}

\section{Introduction}
In this article, we give a new geometric characterization of controllability for a generalized model of the nonholonomic integrator and study some 
minimum energy optimal control problems on these models. 
The relationship between optimal control problems for example, minimum energy problems on nonlinear systems and geometric problems on Riemannian manifolds (such as geodesics)
is well known in the control literature due to seminal works of Brockett (\cite{broc2} and the references therein). The celebrated prototype of a 
nonlinear control system to understand these connections is the nonholonomic integrator or the Brockett integrator. In this article, 
we explore further into this prototypical example and its variants and show how some of the optimal control problems are analogous to 
classical electrodynamics problems such as force acting on a particle in an electromagnetic field. The 
optimal state transfer of the nonholonomic integrator and general nonholonomic systems using sinusoids was demonstrated in the works of Murray and 
Sastry \cite{murray, sastry} with applications in motion planning. Motion planning has since been an active research area as can be seen in the 
works of \cite{lavalle,bel,liu,liu1,henninger,biggs,yang,car,zuy} and the references therein.  
We used orthogonal polynomials such as Legendre and Chebyshev polynomials for steering the nonholonomic 
integrator in \cite{steering} and showed that these orthogonal polynomials can serve as optimal inputs for appropriate cost functions using 
Sturm-Liouville theory. 

	The following nonlinear control system 
	\begin{equation}
 \dot{x}_1=u_1,\; \dot{x}_2=u_2,\;\dot{x}_3=-x_2u_1+x_1u_2\label{areaform}
\end{equation}
is known as the nonholonomic integrator. Notice that the dynamics in the third state co-ordinate is actually a differential $1-$form 
$dx_3=x_1dx_2-x_2dx_1$ (\cite{broc1}). In \cite{broc1}, nonlinear systems with more general differential $1-$forms in $x_1,x_2$ co-ordinates 
are considered. In specific, \cite{broc1} analyzes optimal control problem on systems of dimensions greater than three, where the dynamics on the 
first two state components were defined as in $(\ref{areaform})$ and dynamics on the remaining state components were defined using different 
$1-$forms. These type of systems arise in robotics and motion planning (\cite{murrayli, gurv,yang,car,zuy} and the references therein.  
The following more general form  
\begin{equation}
 \dot{x}_1=u_1,\; \dot{x}_2=u_2,\;\dot{x}_3=f_1(x_1,x_2)u_1+f_2(x_1,x_2)u_2\label{genareaform}
\end{equation}
was considered in \cite{gurv,yang,car,zuy}. Borrowing these ideas, we consider $(\ref{genareaform})$ and its various generalizations in this article. These 
models are important because they provide a canonical form a wider class of nonholonomic control systems and the more general nonholonomic systems can be 
better understood by studying these specific nonholonomic systems. 
%
%
%

We use notions such as the curl of a vector field from multi-variable calculus to give necessary and sufficient conditions for controllability of 
$(\ref{genareaform})$ and its generalizations. We show that minimum energy optimal control problems for these models can be identified with the 
classical electrodynamics problem of a particle in an electromagnetic field. We then give another characterization of controllability using holomorphic functions 
from complex analysis. \\
{\bf Organization}: In the next section, we give some preliminaries to be used in this paper. In Section $\ref{sec:ctrl}$, we obtain a necessary 
and sufficient condition of controllability for generalizations of the nonholonomic integrator using the curl operator. Then, in Section $\ref{sec:electro}$, we 
explore the relationship between some minimum energy optimal control problems on the nonholonomic integrator and classical electrodynamics. In 
Section $\ref{sec:complexanal}$, we study controllability of the general nonholonomic integrator using tools from complex analysis. \\
{\bf Notation}: The scalars and scalar valued functions are denoted by small-face letters, vectors and vector valued functions are denoted by 
bold-face letters and matrices and matrix valued functions are denoted by capital letters. The gradient operator on a scalar valued function $\phi$ is 
denoted as $\nabla \phi$, the curl operator on a vector field $\bold{f}$ in $\mathbb{R}^3$ or $\mathbb{R}^2$ is denoted by $\nabla \times \bold{f}$ 
and the divergence operator is denoted by $\nabla.\bold{f}$. The closed loop integral over a closed curve $\gamma$ is denoted by $\oint_{\gamma}$ and 
the surface integral over a surface $S$ is denoted by $\int\int_{S}$. The line element on a manifold is denoted by $ds$. 

\section{Preliminaries}\label{prelim}
 We refer to Equation $(\ref{areaform})$ as the nonholonomic integrator on $\mathbb{R}^2$ for reasons which will become clear later. Notice that if 
$x_1(0)=x_1(1)$ and $x_2(0)=x_2(1)$, then the variable $x_3$ measures the area formed the projection of the state trajectory on $x_1-x_2$ plane 
(follows from Green's theorem). This 
gives some idea why we refer to this system as the nonholonomic integrator on $\mathbb{R}^2$. 
We refer to model considered in Equation $(\ref{genareaform})$ 
as the general nonholonomic integrator or the general nonholonomic integrator on $\mathbb{R}^2$ associated with a vector field $\bold{f}=(f_1,f_2)$ 
on $\mathbb{R}^2$. 


The following system is refereed as the generalized nonholonomic integrator on $\mathbb{R}^m$.
\begin{eqnarray}
 \dot{x}_i&=&u_i,\; i=1,\ldots,m,\nonumber\\
 \dot{x}_{ij}&=& x_iu_j-x_ju_i,\; i<j=1,\ldots,m.\label{gennonhol} 
\end{eqnarray}
Suppose $x_i(0)=x_i(T)=0$, $\forall i=1,\ldots,m$. Then, the co-ordinates ${x}_{ij}$ measure the area of the closed curve obtained by the 
projection of the state trajectory onto $x_i-x_j$ plane. 

We also consider the following form of the nonholonomic integrator on $\mathbb{R}^3$ in the sequel
\begin{eqnarray}
 &&\dot{x}_1=u_1,\; \dot{x}_2=u_2,\;\dot{x}_3=u_3,\nonumber\\
 &&\dot{x}_4=f_1(x_1,x_2,x_3)u_1+f_2(x_1,x_2,x_3)u_2+f_3(x_1,x_2,x_3)u_3.\label{genareaform1}
\end{eqnarray}
Moreover, we also consider the following generalization of $(\ref{gennonhol})$ 
\begin{eqnarray}
 \dot{x}_i&=&u_i,\; i=1,\ldots,m,\nonumber\\
 \dot{x}_{ij}&=& f_i(x_i,x_j)u_i+f_j(x_i,x_j)u_j,\; i<j=1,\ldots,m.\label{gennonhol1} 
\end{eqnarray}
We refer the reader to \cite{broc1} and \cite{sastry} for more generalizations of the nonholonomic integrator. 
As far as steering of the classical nonholonomic integrator and generalized nonholonomic integrator is concerned, \cite{murray} gave a steering 
algorithm using sinusoids which also holds for a wider class of nonholonomic systems such as the ones defined above. 

We briefly mention the following example from \cite{sastry} which gives optimal sinusoidal inputs for the nonholomorphic integrator when the cost function is the 
minimum input energy function. 
 \begin{example}[\cite{sastry}]\label{sasex}
 For the system defined by $(\ref{areaform})$, we want to find the minimum energy input to drive the state from the origin to a specified point $(0,0,a)$ from $t=0$ to $t=1$. The cost function 
 is $J=\int_0^1(u_1^2+u_2^2)dt$ subject to the system dynamics. 
 Using system equations to eliminate $u_1$ and $u_2$, we obtain the cost function $\int_0^1(\dot{x}_1^2+\dot{x}_2^2)dt$ subject to 
 $\dot{x}_3-x_1\dot{x}_2+x_2\dot{x}_1=0$. Therefore, 
 the augmented cost function is 
 \begin{eqnarray}
  J_a=\int_0^1(\dot{x}_1^2+\dot{x}_2^2+p(t)(\dot{x}_3-x_1\dot{x}_2+x_2\dot{x}_1))dt.\nonumber
 \end{eqnarray}
Applying the first order necessary conditions from calculus of variations, 
we obtain 
$p(t)=c$ and 
\begin{eqnarray}
 \ddot{x}_1+c\dot{x}_2&=&0\nonumber\\
 \ddot{x}_2-c\dot{x}_1&=&0.\nonumber
\end{eqnarray}
Now using $\dot{x}_1=u_1$ and $\dot{x}_2=u_2$, we have the following first order ode
\begin{eqnarray}
 \dot{\left[ \begin{array}{c} {u}_1\\u_2\end{array} \right]}&=&\left[ \begin{array}{cc} 0&-c\\c&0\end{array} \right]
 \left[ \begin{array}{c} {u}_1\\u_2\end{array} \right]\nonumber\\
 \Rightarrow \left[ \begin{array}{c} {u}_1\\u_2\end{array} \right]&=&\left[ \begin{array}{cc} \cos ct&-\sin ct\\\sin ct& \cos ct\end{array} \right]
 \left[ \begin{array}{c} {u}_1(0)\\u_2(0)\end{array} \right].\nonumber
\end{eqnarray}
We need to find $\bold{u}(0)$ and $c$ using initial and final conditions. Let's write $\dot{\bold{u}}=H\bold{u}$ for first order equations in $u_1,u_2$. 
Hence, $\bold{u}(t)=e^{Ht}\bold{u}(0)$. Note that $e^{Ht}$ is orthogonal, hence, the norm of $\|\bold{u}(t)\|=\|\bold{u}(0)\|$ remains constant for 
all time. 
 From the terminal conditions, it follows that $c=2n\pi$ where $n=0,\pm 1,\pm 2,\ldots$. 
Suppose $a>0$, then the cost is 
minimum when $n=1$ and $\|\bold{u}\|=2\pi a$ 
with the direction of $\bold{u}$ being arbitrary.

For an arbitrary terminal time $T$, 
it turns out that $cT=2n \pi$. 
Thus, for $n=1$, $c=\frac{2\pi}{T}$ and  
\begin{equation}
 \left[ \begin{array}{c} {u}_1(t)\\u_2(t)\end{array} \right]=\left[ \begin{array}{cc} \cos \frac{2\pi}{T} t&-\sin \frac{2\pi}{T} t\\\sin \frac{2\pi}{T} t& 
 \cos \frac{2\pi}{T} t\end{array} \right]
 \left[ \begin{array}{c} {u}_1(0)\\u_2(0)\end{array} \right].\nonumber
\end{equation}
Let $u_i(0)=\sqrt{\frac{ca}{2}}$, $i=1,2$. 
Therefore, with sinusoidal inputs of appropriate frequencies, one can always 
steer the system from the origin to any point $(0,0,a)$ in time $T$. The frequencies are chosen depending upon the terminal time $T$ so that for $x_1$ and 
$x_2$, we are integrate the sinusoids over the full period. 

 \end{example}
{\bf Holomorphic functions and Cauchy's integral formula}: A function $F:\mathbb{C}\rightarrow \mathbb{C}$ is called holomorphic if it is complex differentiable 
at each point in $\mathbb{C}$. Let $z=x_1+ix_2\in\mathbb{C}$ and $F(z)=F_1(x_1,x_2)+iF_2(x_1,x_2)$. Then, for a holomorphic function $F$, its real and 
imaginary parts $F_1,F_2$ satisfy Cauchy-Riemann equations given by (\cite{strang})
\begin{eqnarray}
 \frac{\partial F_1}{\partial x_1}=\frac{\partial F_2}{\partial x_2},\;\frac{\partial F_2}{\partial x_1}=-\frac{\partial F_1}{\partial x_2}.\label{CReqn}
\end{eqnarray}

Let $U\subseteq \mathbb{C}$ be an open subset and $F$ be a holomorphic function on $U$. Let $\gamma\subset U$ be a closed 
curve. Then, Cauchy's integral theorem says that $\oint_{\gamma}F(z)dz=0$. 
Let $a$ be a point in the interior of the curve $\gamma$. Then, Cauchy's integral formula says that (\cite{strang})
\begin{equation}
 F(a)=\frac{1}{2\pi i}\oint_{\gamma}\frac{F(z)}{z-a}dz\label{cif}
\end{equation}
which can be proved using Cauchy's integral theorem. 

\section{Characterization of controllability using the curl operator} \label{sec:ctrl}
Consider the system $(\ref{genareaform})$.  
By Green's theorem, $x_3$ measures $\int\int(\frac{\partial f_2}{\partial x_1}-\frac{\partial f_1}{\partial x_2})dx_1dx_2$ over the area enclosed by the loop 
obtained by the projection of the state trajectory on $\mathbb{R}^2$. We can measure the divergence or the curl of a vector field using the $x_3$ coordinate. 
In specific, one can define a vector field on $\mathbb{R}^2$. 
Suppose the projection of the state trajectory on $\mathbb{R}^2$ forms a loop. The $x_3$ coordinate measures 
the curl of the vector field. We now show how it is related to controllability. 

The following theorem gives necessary and sufficient conditions for controllability of $(\ref{genareaform1})$ in terms of the geometry of the underlying 
vector field $\bold{f}$. Notice that one only needs to check arbitrary state transfer of the state variable $x_4$. 
\begin{theorem}\label{geomctrlb}
 Consider system $(\ref{genareaform1})$ and let $\bold{f}=(f_1,f_2,f_3)$ be a continuously differentiable vector field on $\mathbb{R}^3$. 
 The following are equivalent
 \begin{enumerate}
  \item The system $(\ref{genareaform1})$ is controllable.
  \item There exists a closed loop $\gamma \in \mathbb{R}^3$ such that 
  the line integral $\oint_{\gamma}\bold{f}.d\bold{x}\neq0$. 
  \item $\nabla \times \bold{f}\neq0$ on $\mathbb{R}^3$. 
 \end{enumerate}
\end{theorem}
\begin{IEEEproof}
 $(1)\Rightarrow (2)$ 
 Suppose $\oint_{\gamma}\bold{f}.d\bold{x}=0$ for every closed loop $\gamma \in \mathbb{R}^3$, then one cannot do a state transfer from the origin to $(0,0,0,a)$ hence, the system is uncontrollable. 
 $(2)\Rightarrow (3)$ follows from Stokes' theorem. 
 Suppose $(3)$ is satisfied. Let $S$ be a some two dimensional surface in $\mathbb{R}^3$ with the boundary $\gamma$ such that $\nabla \times \bold{f}\neq0$ 
 on $S$ and the surface integral $\iint_{S}(\nabla \times \bold{f}).dS\neq0$. 
 Since $x_1,x_2,x_3$ are controllable, one can choose $u_i$ $(i=1,2,3)$ such that projection of the 
 state trajectory on $\mathbb{R}^3$ is given by $\gamma$. Now since $\iint_{S}(\nabla \times \bold{f}).dS\neq0$, using Stokes' theorem, $x_4$ 
 can also be steered which proves controllability.
\end{IEEEproof}
\begin{corollary}
Consider system $(\ref{genareaform})$ and let $\bold{f}=(f_1,f_2)$ be a continuously differentiable vector field on $\mathbb{R}^2$.
 The following are equivalent
 \begin{enumerate}
  \item The system $(\ref{genareaform})$ is controllable.
  \item There exists a closed loop $\gamma \in \mathbb{R}^2$ such that 
  the line integral $\oint_{\gamma}\bold{f}.d\bold{x}\neq0$. 
  \item $\nabla \times \bold{f}\neq0$ on $\mathbb{R}^2$. 
 \end{enumerate}
\end{corollary}
\begin{IEEEproof}
 Follows from the proof of Theorem $\ref{geomctrlb}.$
\end{IEEEproof}
\begin{example}\label{rigid_osc}
 Consider the following system 
 \begin{equation}
  \dot{x}_1=u_1,\;\dot{x}_2=u_2,\;\dot{x}_3=x_2^2u_1-x_1^2u_2\label{rigosc}
 \end{equation}
from \cite{yang,car,zuy} which describes the motion of a planar rigid body with two oscillators. Clearly since $\nabla\times \bold{f}\neq0$, the system 
is controllable. We now give an explicit steering of the system from $(0,0,0)$ at $t=0$ to $(0,0,a)$ at $t = 1$, where $a>0$. 
Suppose $u_1=c_1\cos (2\pi t)$ and $u_2=c_2\sin (2\pi t)$. Therefore, $x_1(t)=\frac{c_1}{2\pi}(\sin (2\pi t))$ and 
$x_2(t)=\frac{c_2}{2\pi}(1-\cos (2\pi t))$. Now $x_3(t)$ is given by
\begin{gather}
x_3(t) = \int_{0}^{t}\frac{c_1c_2^2}{4\pi^2}(1-\cos (2\pi t))^2cos(2\pi t) - \frac{c_1^2c_2}{4\pi^2}(\sin^2(2\pi t))\sin (2\pi t)dt \\
\Rightarrow x_3(1) = \int_{0}^{1}\frac{c_1c_2^2}{4\pi^2}(1-\cos (2\pi t))^2cos(2\pi t) - \frac{c_1^2c_2}{4\pi^2}(\sin^3(2\pi t))dt \\
\Rightarrow x_3(1) = a = -\frac{c_1c_2^2}{4\pi^2}.
\end{gather}
Thus, for appropriate choices of $c_1,c_2$, we can steer the system from $(0,0,0)$ to $(0,0,a)$.
\end{example} 
\begin{remark}
Notice that if $\nabla \times \bold{f}=0$, then the system is uncontrollable. Thus, if $\bold{f}$ is a gradient vector field i.e. $\bold{f}=\nabla \phi$ 
for some potential function $\phi$, then $(\ref{genareaform})$ 
and $(\ref{genareaform1})$ are uncontrollable. 

Consider the system defined by $(\ref{genareaform})$. Let $\bold{f}=\left[ \begin{array}{cc} 1&0\\0&-1\end{array} \right]\nabla \phi$. Then, 
$\nabla \times \bold{f}\neq0$ in general. Therefore, one can construct controllable systems using a scalar potential function. Similarly, for systems defined by 
$(\ref{genareaform1})$, we can construct controllable systems using $\bold{f}=H\nabla \phi$ where 
\begin{eqnarray}
 H=\left[ \begin{array}{ccc} 1&0&0\\0&-1&0\\0&0&-1\end{array} \right]\mbox{ or } H=\left[ \begin{array}{ccc} 1&0&0\\0&1&0\\0&0&-1\end{array} \right]\nonumber
\end{eqnarray}
and so on. There are other choices of $H$ as well which ensure that $\nabla\times (H\nabla\phi)\neq0$ apart from the ones given above. 
\end{remark}
\begin{corollary}
Consider the system defined by $(\ref{gennonhol1})$ and for $1\le i < j\le n$, 
let $\bold{F}_{ij}=(f_i(x_i,x_j),f_j(x_i,x_j))$ be continuously differentiable vector fields on $\mathbb{R}^2$. 
 The following are equivalent
 \begin{enumerate}
  \item The system $(\ref{gennonhol1})$ is controllable.
  \item For each $1\le i < j\le m$, there exists a closed loop $\gamma \in \mathbb{R}^2$ such that 
  the line integral $\oint_{\gamma}\bold{F}_{ij}.d\bold{x}\neq0$. 
  \item $\nabla \times \bold{F}_{ij}\neq0$ on $\mathbb{R}^2$. 
 \end{enumerate}

\end{corollary}
\begin{IEEEproof}
 Follows from the proof of Theorem $\ref{geomctrlb}$ and the previous corollary.
\end{IEEEproof}

\begin{remark}
 Notice that since $\mathbb{R}^3$ and $\mathbb{R}^2$ are simply connected, $\nabla \times \bold{f}=0$ $\Leftrightarrow \bold{f}=\nabla \phi$ for some 
 scalar function $\phi$. Therefore, $(\ref{genareaform1})$ is uncontrollable $\Leftrightarrow$ $\bold{f}$ is a gradient vector field.  
\end{remark}
\begin{example}\label{uncex}
 Consider a system 
 \begin{equation}
  \dot{x}_1=u_1,\;\dot{x}_2=u_2,\;\dot{x}_3=\frac{x_1}{x_1^2+x_2^2}u_1+\frac{x_2}{x_1^2+x_2^2}u_2 \nonumber
 \end{equation}
defined over $\mathbb{R}^{3}\setminus{\{0,0,x_3\}}$. Notice that $\nabla \times \bold{f}=0$ over $\mathbb{R}^{3}\setminus{\{0,0,x_3\}}$ and the system is  
uncontrollable. Since the state space is not simply connected, although $\nabla\times \bold{f}=0$, $\bold{f}$ is not a gradient vector field. 
\end{example}
%

Suppose we want to steer $(\ref{genareaform})$ from the origin to $(0,0,a)$. The $x_3$ coordinate is given by 
\begin{eqnarray}
x_3(t)=\int_{0}^{t} (f_1u_1 + f_2u_2)dt.\nonumber
\end{eqnarray}
Choose $u_1,u_2$ as orthogonal polynomials so that $x_1(1)=x_2(1)=0$. 
\begin{eqnarray}
x_3(1)=\int_{0}^{1} (f_1u_1 + f_2u_2)dt.\nonumber
\end{eqnarray}
We need to choose orthogonal $u_1,u_2$ such that the above integral is nonzero. 
To steer from the origin to $(a,b,c)$, choose constant inputs to steer the state to some point say $(a,b,d)$. Then use orthogonal polynomials 
to steer along $x_3$ without affecting $x_1,x_2$. One can similarly steer the nonholonomic integrator $(\ref{genareaform1})$ on $\mathbb{R}^3$. 
(Notice that in the case above, choosing $u_1=-f_2, u_2=f_1$, the motion can be constrained to $x_1-x_2$ plane.) 

%
\begin{example}
	Consider the system given by $(\ref{genareaform})$ where $f_1(x_1,x_2)=x_1^2-x_2^2$ and $f_2(x_1,x_2)=2x_1x_2$. It follows that $\nabla \times \bold{f}=4x_2\neq0$. Therefore, 
	the system is controllable by Theorem $\ref{geomctrlb}$. Suppose we want to steer the system from the origin at $t=0$ to $(0,0,a)$ at $t=1$ where $a>0$. 
	Suppose $u_1=c_1\cos (2\pi t)$ and $u_2=c_2\sin (2\pi t)$. Therefore, $x_1(t)=\frac{c_1}{2\pi}(\sin (2\pi t))$ and 
	$x_2(t)=\frac{c_2}{2\pi}(1-\cos (2\pi t))$. Now $x_3(t)$ is given by
	\begin{gather*}
	x_3(t) = \int_{0}^{t} \bigg(\frac{c_1^2}{4\pi^2}\sin^2(2\pi t) - \frac{c_2^2}{4\pi^2}(1-\cos (2\pi t))^2\bigg) c_1 \cos(2\pi t) + \\
	2\frac{c_1}{2\pi}(\sin (2\pi t))\frac{c_2}{2\pi}(1-\cos (2\pi t))c_2 \sin(2\pi t) dt \\
	\Rightarrow x_3(1) = \int_{0}^{1} \bigg(\frac{c_1^2}{4\pi^2}sin^2(2\pi t) - \frac{c_2^2}{4\pi^2}(1-\cos (2\pi t))^2\bigg) c_1 \cos(2\pi t) +\\
	\frac{c_1c_2^2}{2\pi^2}(\sin^2(2\pi t))(1-\cos (2\pi t)) dt \\
	\Rightarrow x_3(1) = \int_{0}^{1} -\frac{c_1c_2^2}{4\pi^2}(1-\cos (2\pi t))^2 \cos(2\pi t) + \frac{c_1c_2^2}{2\pi^2}(\sin^2(2\pi t)) dt \\
	\Rightarrow x_3(1) = a = \frac{c_1c_2^2}{2\pi^2}.    
	\end{gather*}
	Thus, for appropriate choices of $c_1,c_2$, we can steer the system from $(0,0,0)$ to $(0,0,a)$.
\end{example}
\section{Optimal control on the general nonholonomic integrator and classical electrodynamics}\label{sec:electro}
Consider the minimum energy control problem for the system $(\ref{genareaform})$. Applying Euler-Lagrange equations on the augmented Lagrangian 
$L=\dot{x}_1^2+\dot{x}_2^2+\lambda(\dot{x}_3-f_1\dot{x}_1-f_2\dot{x}_2)$, 
\begin{eqnarray}
 \frac{d}{dt}(2\dot{x}_1-\lambda f_1)=-\lambda(\frac{\partial f_1}{\partial x_1}\dot{x}_1+\frac{\partial f_2}{\partial x_1}\dot{x}_2)
 \Rightarrow 2\ddot{x}_1-\lambda\frac{\partial f_1}{\partial x_2}\dot{x}_2+\lambda\frac{\partial f_2}{\partial x_1}\dot{x}_2=0\nonumber\\
 \frac{d}{dt}(2\dot{x}_2-\lambda f_2)=-\lambda(\frac{\partial f_1}{\partial x_2}\dot{x}_1+\frac{\partial f_2}{\partial x_2}\dot{x}_2)
 \Rightarrow 2\ddot{x}_2-\lambda\frac{\partial f_2}{\partial x_1}\dot{x}_1+\lambda\frac{\partial f_1}{\partial x_2}\dot{x}_1=0.\nonumber
\end{eqnarray}
Therefore, (substituting $\lambda$ for $\lambda/2$) 
\begin{equation}
 \dot{\left[ \begin{array}{c} {u}_1\\u_2\end{array} \right]}=\lambda\left[ \begin{array}{cc} 0&\frac{\partial f_1}{\partial x_2}-\frac{\partial f_2}{\partial x_1}
 \\-\frac{\partial f_1}{\partial x_2}+\frac{\partial f_2}{\partial x_1}&0\end{array} \right]
 \left[ \begin{array}{c} {u}_1\\u_2\end{array} \right].\label{optipcousn_r2}
\end{equation}
Notice that $-\frac{\partial f_1}{\partial x_2}+\frac{\partial f_2}{\partial x_1}$ gives the curl of the vector field $(f_1,f_2)$ on $\mathbb{R}^2$. 
One can choose appropriate form for $\dot{x}_3$ so that one obtains the divergence instead. 
\begin{remark}\label{magfldr2}
 It follows from $(\ref{optipcousn_r2})$ and results from \cite{murray} that 
 if the curl of the vector field $\bold{f}=(f_1,f_2)$ on $\mathbb{R}^2$ is constant, then the optimal inputs for $(\ref{genareaform})$ 
 are given by sinusoids for a state transfer 
 from the origin to $(0,0,a)$.
\end{remark}

Rewriting $(\ref{optipcousn_r2})$ as 
\begin{eqnarray}
  \dot{\left[ \begin{array}{c} {u}_1\\u_2\end{array} \right]}&=&\lambda
  (-\frac{\partial f_1}{\partial x_2}+\frac{\partial f_2}{\partial x_1})\left[ \begin{array}{cc} 0&-1\\1&0\end{array} \right]
 \left[ \begin{array}{c} {u}_1\\u_2\end{array} \right]\nonumber\\
 &=&\lambda
  (-\frac{\partial f_1}{\partial x_2}+\frac{\partial f_2}{\partial x_1})\left[ \begin{array}{c} -{u}_2\\u_1\end{array} \right]\nonumber\\
 &=&\lambda
  (\nabla \times \bold{f})\times \bold{u}\label{magfldeqn1}
\end{eqnarray}
where $\bold{f}=(f_1,f_2)=f_1\hat{i}+f_2\hat{j}$ in the vectorial notation and $\bold{u}=(u_1,u_2)=u_1\hat{i}+u_2\hat{j}$. Notice that 
$(\nabla \times \bold{f})= (-\frac{\partial f_1}{\partial x_2}+\frac{\partial f_2}{\partial x_1})\hat{k}$ (where $\hat{i},\hat{j}$ and $\hat{k}$ are 
unit vectors along $x_1,x_2$ and $x_3$ respectively). Now substituting $\dot{x}_i=u_i$ $(i=1,2)$, one obtains 
\begin{equation}
 \ddot{\bold{x}}=\lambda  (\nabla \times \bold{f})\times \dot{\bold{x}}\label{magfldeqn2}
\end{equation}
where $\bold{x}=(x_1,x_2)$. Observe that $(\nabla \times \bold{f})$ can be thought of as a magnetic field and $\dot{\bold{x}}$ is the velocity. Thus, 
$(\ref{optipcousn_r2})$ can be interpreted as the force acting on a particle in a magnetic field. This relates the optimal control problem to the classical 
electrodynamics. Magnetic field in classical electrodynamics is given by the curl of a vector potential. Any vector field can be decomposed by the 
Helmholtz decomposition into irrotational (curl-free) and solenoidal (divergence free) vector field. Thus, adding $\nabla \phi$ to $\bold{f}$ still 
gives the same dynamics on $x_3$. In other words, only the solenoidal component of the vector field $\bold{f}$ plays a role in solving the optimal 
control problem. For gradient vector fields, $\bold{f}=\nabla \phi$, there is no motion possible. For details on classical mechanics, we refer the reader 
to \cite{arnold} and for classical electrodynamics, we refer the reader to \cite{griffiths}. 



Consider the minimum energy control problem of minimizing $\int_0^1(u_1^2+u_2^2+u_3^2)dt$ on $(\ref{genareaform1})$. The augmented Lagrangian is 
$L=\dot{x}_1^2+\dot{x}_2^2+\dot{x}_3^2+\lambda(\dot{x}_4-f_1\dot{x}_1-f_2\dot{x}_2-f_3\dot{x}_3)$. Using Euler-Lagrange equations, one obtains
\begin{eqnarray}
 2\ddot{x}_1&=&\lambda(\frac{\partial f_1}{\partial x_2}-\frac{\partial f_2}{\partial x_1})\dot{x}_2+
 \lambda(\frac{\partial f_1}{\partial x_3}-\frac{\partial f_3}{\partial x_1})\dot{x}_3\label{x1eq}\\
 2\ddot{x}_2&=&\lambda(\frac{\partial f_2}{\partial x_1}-\frac{\partial f_1}{\partial x_2})\dot{x}_1+
 \lambda(\frac{\partial f_2}{\partial x_3}-\frac{\partial f_3}{\partial x_2})\dot{x}_3\label{x2eq}\\
 2\ddot{x}_3&=&\lambda(\frac{\partial f_3}{\partial x_1}-\frac{\partial f_1}{\partial x_3})\dot{x}_1+
 \lambda(\frac{\partial f_3}{\partial x_2}-\frac{\partial f_2}{\partial x_3})\dot{x}_2\label{x3eq}
\end{eqnarray}
and $\lambda$ is a constant. Substituting $\lambda$ for $\lambda/2$, the above equations can be written in the matrix form as
\begin{eqnarray}
 \ddot{\left[ \begin{array}{c} {x}_1\\x_2\\x_3\end{array} \right]}&=&
 \lambda\left[ \begin{array}{ccc} 0&\frac{\partial f_1}{\partial x_2}-\frac{\partial f_2}{\partial x_1}&\frac{\partial f_1}{\partial x_3}-\frac{\partial f_3}{\partial x_1}
 \\\frac{\partial f_2}{\partial x_1}-\frac{\partial f_1}{\partial x_2}&0&\frac{\partial f_2}{\partial x_3}-\frac{\partial f_3}{\partial x_2}\\
 \frac{\partial f_3}{\partial x_1}-\frac{\partial f_1}{\partial x_3}&\frac{\partial f_3}{\partial x_2}-\frac{\partial f_2}{\partial x_3}&0\end{array} \right]
 \left[ \begin{array}{c} \dot{x}_1\\\dot{x}_2\\\dot{x}_3\end{array} \right]\nonumber\\
 &=&\lambda (\nabla \times \bold{f})\times \dot{\bold{x}}.\label{optipcousn_r3}
\end{eqnarray}
Again, we have an equation of motion of a particle in a magnetic field in $\mathbb{R}^3$. 

\begin{remark}
 Consider the steering problem in $(\ref{genareaform1})$ where one wants to steer from the origin to $(0,0,0,a)$. 
Then,  
\begin{equation}
 x_4(1)=\int_0^1(f_1u_1+f_2u_2+f_3u_3)dt=\oint f_1dx_1+f_2dx_2+f_3dx_3\nonumber
\end{equation}
where the closed loop integral is over the closed curve obtained by projection of the state trajectory in $\mathbb{R}^4$ on to $\mathbb{R}^3$. This 
can be interpreted as the work done by the vector field $\bold{f}$ along the curve. Thus, $x_4$ measures the work done by $\bold{f}$ along the 
projected curve. Furthermore, by Stokes' theorem, $x_4$ also measures the flux of the magnetic field $B=\nabla\times \bold{f}$ passing through any  
surface whose boundary is given by the closed curve obtained above by the projection of the state trajectory from the origin to $(0,0,0,a)$ on $\mathbb{R}^3$. 

\end{remark}

\begin{remark}
 Consider a revised Lagrangian $\hat{L}=\dot{x}_1^2+\dot{x}_2^2+\dot{x}_3^2-\lambda(f_1\dot{x}_1+f_2\dot{x}_2+f_3\dot{x}_3)$ 
where we have used dynamics in $x_4$ to conclude that the Lagrange multiplier $\lambda$ is constant. Notice that for a particle in a magnetic field 
$\bold{B}=\nabla \times \bold{A}$ moving with velocity $\dot{\bold{x}}$, the Lagrangian is given by $\frac{1}{2}m\dot{\bold{x}}.\dot{\bold{x}}+q
\dot{\bold{x}}.\bold{A}$. In the case of the nonholonomic integrator on $\mathbb{R}^3$, $\bold{A}=\bold{f}$. Therefore, the revised Lagrangian 
can be identified with the Lagrangian for electrodynamics. Thus, the optimal control problem is also a classical mechanics problem. Now for the revised 
Lagrangian, the Hamiltonian is preserved and we have Hamiltonian dynamics. 

\end{remark}

\subsection{Incorporating a drift term in the nonholonomic integrator and its relation to the force on a particle in an electromagnetic field}
Consider the following system with a drift 
\begin{eqnarray}
 &&\dot{x}_1=u_1,\; \dot{x}_2=u_2,\;\dot{x}_3=u_3,\nonumber\\
 &&\dot{x}_4=g(x_1,x_2,x_3)+f_1(x_1,x_2,x_3)u_1+f_2(x_1,x_2,x_3)u_2+\nonumber\\
 &&f_3(x_1,x_2,x_3)u_3\label{genareaform2}
\end{eqnarray}
and the minimum energy control problem of minimizing $\int_0^1(u_1^2+u_2^2+u_3^2)dt$. The augmented Lagrangian is 
$L=\dot{x}_1^2+\dot{x}_2^2+\dot{x}_3^2+\lambda(\dot{x}_4-g-f_1\dot{x}_1-f_2\dot{x}_2-f_3\dot{x}_3)$. Using Euler-Lagrange equations, one obtains
\begin{eqnarray}
 2\ddot{x}_1&=&-\lambda\frac{\partial g}{\partial x_1}+\lambda(\frac{\partial f_1}{\partial x_2}-\frac{\partial f_2}{\partial x_1})\dot{x}_2+
 \lambda(\frac{\partial f_1}{\partial x_3}-\frac{\partial f_3}{\partial x_1})\dot{x}_3\label{x1eqd}\\
 2\ddot{x}_2&=&-\lambda\frac{\partial g}{\partial x_2}+\lambda(\frac{\partial f_2}{\partial x_1}-\frac{\partial f_1}{\partial x_2})\dot{x}_1+
 \lambda(\frac{\partial f_2}{\partial x_3}-\frac{\partial f_3}{\partial x_2})\dot{x}_3\label{x2eqd}\\
 2\ddot{x}_3&=&-\lambda\frac{\partial g}{\partial x_3}+\lambda(\frac{\partial f_3}{\partial x_1}-\frac{\partial f_1}{\partial x_3})\dot{x}_1+
 \lambda(\frac{\partial f_3}{\partial x_2}-\frac{\partial f_2}{\partial x_3})\dot{x}_2\label{x3eqd}
\end{eqnarray}
where $\lambda$ is a constant. Now the above equations can be expressed as 
\begin{equation}
 \ddot{\bold{x}}=-\lambda\nabla g + \lambda(\nabla \times \bold{f})\times \dot{\bold{x}}\label{elecdyn}
\end{equation}
which can be identified with the classical force equation $\bold{F}=q\bold{E}+m\bold{v}\times\bold{B}$. This is again a Hamiltonian system where the Lagrangian 
is defined by $L=\dot{x}_1^2+\dot{x}_2^2+\dot{x}_3^2-\lambda(g+f_1\dot{x}_1+f_2\dot{x}_2+f_3\dot{x}_3)$. 
\begin{remark}
 One can show that $(\ref{genareaform2})$ is controllable $\Leftrightarrow$ $\nabla \times \bold{f}\neq0$ using similar arguments used in Theorem  
 $\ref{geomctrlb}$. If $\nabla \times \bold{f}=0$, then $x_4$ cannot be steered arbitrarily. 
\end{remark} 

%
\begin{remark}
%
 Consider a $4-$vector potential $(g,\bold{f})$ in analogy with the $4-$vector potential $(\phi/c,\bold{A})$ in electrodynamics. The Lorenz gauge 
 condition is given by $\frac{1}{c^2}\frac{\partial \phi}{\partial t}+\nabla . \bold{A}=0$. Thus, $\phi/c=g$ and $\bold{A}=\bold{f}$. Since only 
 $\nabla\times\bold{f}$ decides the controllability of the system, one can ignore the curl free part in the Helmholtz decomposition and $\bold{f}$ 
 can be assumed to be solenoidal. Therefore, $\nabla.\bold{f}=0$. Moreover, $\frac{\partial g}{\partial t}=0$ which implies that Lorenz gauge conditions 
 are satisfied in the above case as well. If $g$ has an explicit time dependence, then $\nabla.\bold{f}\neq 0$ but we do not consider this case here. 
%

Note that in electrodynamics, $\bold{B}=\nabla \times \bold{A}$ and $\bold{E}=-\frac{\partial A}{\partial t}-\nabla \phi$. 
Now for the control systems considered above, if $g=0$ and $\bold{f}$ has no explicit time dependence, then $\bold{E}=0$ and 
$\bold{B}=\nabla\times \bold{f}$. When $g\neq0$ or if $\bold{f}$ is time dependent, then $\bold{E}\neq0$. Notice that controllability of 
$(\ref{genareaform})$, $(\ref{genareaform1})$ and $(\ref{genareaform2})$ can be related to 
the presence of a magnetic field. In the absence of magnetic fields i.e., $\nabla \times \bold{f}=0$, the system becomes uncontrollable. 
\end{remark}
\subsection{State dependent cost function and electrodynamics analogy} 
Consider the following optimal control problem on $(\ref{genareaform1})$ involving a quadratic cost on the first three components of the state 
\begin{eqnarray}
J&=&\int_{0}^{1}(x_1^2+x_2^2+x_3^2+u_1^2+u_2^2+u_3^2)dt\nonumber\\
&=&\int_{0}^{1}(x_1^2+x_2^2+x_3^2+\dot{x}_1^2+\dot{x}_2^2+\dot{x}_3^2)dt.\label{sdcost} 
\end{eqnarray}
The augmented Lagrangian is 
$L=x_1^2+x_2^2+x_3^2+\dot{x}_1^2+\dot{x}_2^2+\dot{x}_3^2+\lambda(\dot{x}_4-f_1\dot{x}_1-f_2\dot{x}_2-f_3\dot{x}_3)$. Using Euler-Lagrange equations, one obtains
\begin{eqnarray}
 2\ddot{\bold{x}}=2\bold{x}+\lambda(\nabla\times\bold{f})\times \dot{\bold{x}}\nonumber
\end{eqnarray}
where $\bold{x}\in\mathbb{R}^3$ and $\lambda \in \mathbb{R}$. One could consider a state dependent term $g(x_1,x_2,x_3)>0$ instead of the quadratic term 
$x_1^2+x_2^2+x_3^2$ in the cost function $(\ref{sdcost})$ 
to obtain 
\begin{eqnarray}
 2\ddot{\bold{x}}=\nabla g(\bold{x})+\lambda(\nabla\times\bold{f})\times \dot{\bold{x}}.\nonumber
\end{eqnarray}
Thus, the electrodynamics analogy can be obtained for the general nonholonomic integrator with a drift term or for the general nonholonomic integrator with 
a state dependent term $g(x_1,x_2,x_3)>0$ in the cost function as shown in Equation $(\ref{sdcost})$ where $g(x_1,x_2,x_3)=x_1^2+x_2^2+x_3^2$. It is clear 
that all these optimal control problems on general versions of nonholonomic integrator can be identified with solving an electrodynamics problem. 
\subsection{Optimal control of a planar rigid body with two oscillators}
We now demonstrate with an example, how to reduce the minimum energy optimal control problem (which is also related to finding the trajectory 
of a particle in a magnetic field) to solving an elliptic integral.  

\begin{example}
Consider the following system which describes the motion of a planar rigid body with two oscillators. (Refer Example (\ref{rigid_osc})).
\begin{equation}
\dot{x}_1=u_1,\;\dot{x}_2=u_2,\;\dot{x}_3=x_2^2u_1-x_1^2u_2.
\end{equation}
Here, we want to optimize the following cost function given by 
\begin{gather}
J = \int_{0}^{1} (u_1^2 + u_2^2) dt
\end{gather}
subject to the constraints that $(x_1(0),x_2(0),x_3(0))=(0,0,0)$ and $(x_1(1),x_2(1),x_3(1))=(0,0,c)$. 
Here, $\mathbf{f} = (x_2^2, -x_1^2)$ and the equations of motion are 
\begin{gather}
2\ddot{\mathbf{x}} = \lambda (\nabla \times \mathbf{f})\times \dot{x} \\
\ddot{x_1} = \lambda(x_1+x_2)\dot{x_2} \\
\ddot{x_2} = -\lambda(x_1+x_2)\dot{x_1}
\end{gather}
Now, we define a change variables as follows $y = x_1-x_2$, $z = x_1 +x_2$, then equations of motion reduce to 
\begin{gather}
\ddot{y} = \lambda z\dot{z} \label{eqex1}\\
\ddot{z} = -\lambda z\dot{y} \label{eqex2}\\
\dot{y}\ddot{y} + \dot{z}\ddot{z} = 0 \label{eqex3}\\
\dot{y}^2 + \dot{z}^2 = r^2 \label{eqex4}
\end{gather}
where r is a constant. Thus, integrating Equation $(\ref{eqex1})$ and then substituting the result into Equation $(\ref{eqex4})$, we obtain 
\begin{gather}
\dot{y} = \lambda \frac{z^2}{2} + c \\
\dot{z}^2 = r^2 - (\lambda \frac{z^2}{2} + c)^2 \\
dt = \frac{dz}{\sqrt{(r+c+\lambda \frac{z^2}{2})(r-c-\lambda \frac{z^2}{2})}}. \\
{\mbox{Let }} z = \sqrt{\frac{2(r-c)}{\lambda}} \sin \theta, \mbox{ then,} \nonumber \\
dt = \sqrt{\frac{2(r-c)}{\lambda}} \frac{\cos\theta d\theta}{\cos\theta\sqrt{(r-c)(r+c+(r-c)\sin^2\theta)}} \\
dt = \sqrt{\frac{2}{(r+c)\lambda}}\frac{ d\theta}{\sqrt{1+\frac{r-c}{r+c}\sin^2\theta}} .
\end{gather}
Now let $\kappa = \sqrt{\frac{r-c}{r+c}}$, we have,
\begin{gather}
dt = \sqrt{\frac{2}{(r+c)\lambda}}\frac{d\theta}{\sqrt{1+\kappa^2 \sin^2\theta}}. \\
\mbox{Let } u = \frac{\pi}{2} - \theta, \mbox{ then,} \nonumber \\
\sqrt{\frac{(r+c)\lambda}{2}} dt = -\frac{du}{\sqrt{1+\kappa^2(1-\sin^2u)}} = -\frac{du}{\sqrt{1+\kappa^2}\sqrt{1-\frac{\kappa^2}{\kappa^2+1}\sin^2u}} \\
 -\sqrt{\frac{(r+c)\lambda(1+\kappa^2)}{2}} dt = \frac{du}{\sqrt{1+\kappa^2(1-\sin^2u)}} = \frac{du}{\sqrt{1-\frac{\kappa^2}{\kappa^2+1}\sin^2u}} \\
-\sqrt{\frac{(r+c)\lambda(1+\kappa^2)}{2}}t + b = F\bigg(u|\frac{\kappa^2}{\kappa^2+1}\bigg) \label{jacobielliptic}
\end{gather}
where, in the Equation (\ref{jacobielliptic}), $F(\psi|\kappa^2)$ represents the incomplete elliptic integral of the first kind. 
Thus, the optimal control problem reduces to solving the elliptic integral.  
For further 
details on elliptic integrals, we refer the reader to 
\cite{byrd2013handbook}.
  
\end{example}
\section{Complex analytic characterization of controllability}\label{sec:complexanal}
Consider the complex function $F(x_1+ix_2)=f_2(x_1,x_2)+if_1(x_1,x_2)$ where $f_i(x_1,x_2)$ are both real valued functions corresponding to the vector field 
$\bold{f}=(f_1,f_2)$. 
Let $z=x_1+ix_2$ and $u_{\mathbb{C}}=u_1+iu_2$. Let $\gamma$ be a closed curve in $\mathbb{C}$. Then, 
\begin{eqnarray}
 \oint F.u_{\mathbb{C}} dt&=&\oint ((f_2u_1-f_1u_2)+i(f_1u_1+f_2u_2)) dt\nonumber\\
 &=&\oint ((f_2dx_1-f_1dx_2)+i(f_1dx_1+f_2dx_2)).\label{cplxanalexp}
\end{eqnarray}
The function $F$ is defined in such a way that the imaginary part of $F.u_{\mathbb{C}}$ can be identified with the dynamics in the $x_3$ variable associated 
with the control system $(\ref{genareaform})$.
Thus, applying Green's theorem, the real part of the line integral $\oint F.u_{\mathbb{C}} dt$ can be identified with the divergence of $\bold{f}$ 
(since $\oint_{\gamma} (f_2dx_1-f_1dx_2)= \int\int_S(-\nabla.\bold{f})dS$, $S$ being the area enclosed by the closed curve $\gamma$) whereas; 
the imaginary part of the line integral can be identified with the curl of $\bold{f}$ 
(since $\oint_{\gamma} (f_1dx_1+f_2dx_2)=\int\int_S(\nabla \times \bold{f})dS$).  
\begin{lemma}
 Consider $(\ref{genareaform})$ and let $F(x_1+ix_2)=f_2(x_1,x_2)+if_1(x_1,x_2)$. If $F$ is holomorphic, then $(\ref{genareaform})$ is uncontrollable. 
\end{lemma}
\begin{IEEEproof}
 The proof follows from Cauchy's integral theorem since for holomorphic functions, the integral over a closed loop in the complex plane is zero. Thus, the 
 $x_3$ co-ordinate is uncontrollable.
\end{IEEEproof}

\begin{lemma}
 Consider $(\ref{genareaform})$ and let $F(x_1+ix_2)=f_2(x_1,x_2)+if_1(x_1,x_2)$ such that $F$ is not holomorphic. Let $\gamma$ be a closed loop in the complex plane enclosing the 
 origin such that $F$ has a pole in the region enclosed by $\gamma$ and suppose $\gamma$ has a nonzero winding number. If the residue of $F$ at the 
 pole inside $\gamma$ is a nonzero real number, then $(\ref{genareaform})$ is controllable. 
\end{lemma}
\begin{IEEEproof}
 Since $x_1,x_2$ are controllable, one can choose $u_1,u_2$ such that the projection of the state trajectory on the complex plane is given by $\gamma$. 
 Notice that $\dot{x}_3=f_1u_1+f_2u_2$ and by the residue theorem from complex analysis, $x_3$ can be steered if the residue of $F$ at the 
 pole inside $\gamma$ is nonzero real number. 
\end{IEEEproof}
\begin{example}
 Consider the classical nonholonomic integrator with $f_1=-x_2$ and $f_2=x_1$. Therefore, $F=f_2+if_1=x_1-ix_2$. 
 In fact, $F(z)=\bar{z}$ is the complex conjugate of $z$. It can be easily checked that 
 the Cauchy-Riemann equations are not satisfied and the function $F$ is not holomorphic. 
 This agrees with the fact that the classical nonholonomic 
 integrator is controllable. Let $\gamma$ be the unit circle centered at $z=0$. Notice that 
 \begin{eqnarray}
  2\pi i=\oint_\gamma F.u_{\mathbb{C}}=\oint_\gamma (x_1u_1+x_2u_2)+i(x_1u_2-x_2u_1)dt.\nonumber
 \end{eqnarray}
 This implies that the purely imaginary part of 
 $F.u_{\mathbb{C}}$ which captures the dynamics of $x_3$ variable of the nonholonomic integrator is controllable. However,  
 $\oint_\gamma (x_1u_1+x_2u_2)=0$ for every closed curve $\gamma$. Therefore, if $\dot{x}_3=(x_1u_1+x_2u_2) = \text{Re}\{\bar{z}u_{\mathbb{C}}\}$, 
 then the system is not controllable as the closed loop complex integral is always purely imaginary. This can also be verified using the fact that 
 for $\dot{x}_3=x_1u_1+x_2u_2$, $\bold{f}=(x_1,x_2)$ and $\nabla \times \bold{f}=0$ which implies uncontrollability. Thus, uncontrollability in this case is 
 a consequence of the curl of $\bold{f}$ being zero and from complex analytic viewpoint, it follows from the residue theorem.
\end{example}
\begin{example}
 Suppose $F=\frac{1}{z}=\frac{x_1-ix_2}{x_1^2+x_2^2}$ which is not holomorphic and has a pole at the origin. 
 The integral along the closed loop (which is the unit circle)
 using Cauchy's integral formula is given by $2\pi i$. 
 Consider a system 
 \begin{equation}
  \dot{x}_1=u_1,\;\dot{x}_2=u_2,\;\dot{x}_3=\frac{x_1}{x_1^2+x_2^2}u_2-\frac{x_2}{x_1^2+x_2^2}u_1 \nonumber
 \end{equation}
defined over $\mathbb{R}^{3}\setminus{\{0,0,x_3\}}$. Complexify this system using $z=x_1+ix_2$ and $u_{\mathbb{C}}=u_1+iu_2$. Suppose 
 $F=\frac{1}{z}=\frac{x_1-ix_2}{x_1^2+x_2^2}=f_2+if_1$, therefore, $F.u_{\mathbb{C}}=\frac{x_1u_1+x_2u_2}{x_1^2+x_2^2}+i\frac{(-x_2u_1+x_1u_2)}{x_1^2+x_2^2}$. Let $\gamma$ be a simple closed curve in $\mathbb{C}$, which encloses the origin and let $T_p$ be the time period of curve traversed by the chosen inputs $(u_1,u_2)$.  
 Let $\dot{z}=u_{\mathbb{C}}$ and $\dot{x}_3=$ Im$(F.u_{\mathbb{C}})$. It follows that 
 \begin{eqnarray}
  x_3(nT_p)-x_3(0)=\oint_{\gamma}(\frac{x_1}{x_1^2+x_2^2}u_2-\frac{x_2}{x_1^2+x_2^2}u_1)dt=\mbox{Im}\oint_{\gamma}(F.u_{\mathbb{C}})=2n\pi.\nonumber
 \end{eqnarray}
 Therefore, the system is controllable. On the other hand, if $\dot{x}_3=\frac{x_1u_1+x_2u_2}{x_1^2+x_2^2}$, then 
 \begin{eqnarray}
  x_3(nT_p)-x_3(0)=\oint_{\gamma}(\frac{x_1u_1+x_2u_2}{x_1^2+x_2^2})=\mbox{Re}\oint_{\gamma}(F.u_{\mathbb{C}})=0\nonumber
 \end{eqnarray}
and the system is uncontrollable since the real part of the residue is zero and the dynamics in $x_3$ is given by Re$(F.u_{\mathbb{C}})$. 
%
\end{example}
Now consider a nonholonomic system on the complex plane defined as 
\begin{eqnarray}
\dot{z} =  u_{\mathbb{C}} \nonumber\\ 
\dot{w} = F.u_{\mathbb{C}}  \label{fullcomplexform}
\end{eqnarray}
where $z_1=x_1+ix_2$, $u_{\mathbb{C}}=u_1+iu_2$, $F = f_2+if_1$, $w = w_1 + iw_2$ and $x_1,x_2,u_1,u_2,w_1,w_2$ are real variables whereas 
$f_1,f_2$ are functions of real variables $x_1,x_2$. 
In the four dimensional real vector space, this system is represented as follows
\begin{eqnarray}
\dot{x}_1 =  u_1 \\ 
\dot{x}_2 = u_2 \\
\dot{w}_1 = f_2u_1 - f_1u_2  \\
\dot{w}_2 = f_1u_1 + f_2u_2.
\end{eqnarray}

%
%

We now demonstrate how to control a family of complex control systems defined above where the complex function $F=f_2+if_1$ is not holomorphic. 
\begin{example}
	Consider a family of nonholomorphic functions $F(z)=(\bar{z})^{n}$, $n\in\mathbb{N}$ and $n>1$. Consider the nonholonomic system given by 
	\begin{eqnarray}
	\dot{z} =  u_{\mathbb{C}} \nonumber\\ 
	\dot{w} = F.u_{\mathbb{C}} \nonumber
	\end{eqnarray}
	where $z,w$ are defined in $(\ref{fullcomplexform})$. 
	Let $\gamma_a$ be the unit circle centered at $z=a$, then we have (by substituting $v = z-a$ and using the residue theorem)  
	\begin{eqnarray}
	\oint_{\gamma_a} F.u_{\mathbb{C}}=\oint_{\gamma_a} (\bar{z})^{n} dz = \oint_{\gamma_a} (\bar{v} + \bar{a})^n dv = 
	\oint_{\gamma_a} \big(\frac{1}{v} + \bar{a}\big)^n dv = 2 n\pi i (\bar{a}^{n-1}).\nonumber
	\end{eqnarray}
	Thus, the dynamics of $w_2$ can be controlled choosing $a=1$ and appropriate $u_1$ and $u_2$. Note that $w_1$ remains unaffected since the real 
	part of the closed loop integral considered above is zero. Now to control $w_1$, we choose a different complex point say 
	$a = e^{i\frac{\pi}{n}}$. Then, $\bar{a}^{n-1} = e^{-i(n-1)\frac{\pi}{n}} = - e^{i\frac{\pi}{n}} $. Notice that $e^{i\frac{\pi}{n}}$ has both real and 
	imaginary part for all $n \geq 2$, and for $n =2$ it has only the purely imaginary part. Furthermore,  the residues obtained above, when considered as real vectors form 
	a two dimensional real subspace. 
	Thus, the dynamics given by $\dot{w}_1= \text{Re}\{\bar{z}^n.u_{\mathbb{C}}\}$ is controllable and the system with both $w_1$ and $w_2$ is also 
	controllable, as we have two independent real directions associated with the two residues. 
	Now we propose the following algorithm to steer the system from the origin to $(0,0,a,b)$, where $a,b \neq 0$. 
	\begin{enumerate}
		\item Choose the inputs $u_1 = 2\pi \sin(2\pi t)$, $u_2(t) = 2\pi \cos(2\pi t)$, then $x_1 =1-\cos(2\pi t)$, $x_2 = \sin(2\pi t)$. This realizes 
		the curve $\gamma_1$ in $x_1-x_2$ plane and by the discussion above, $w_2(1) = 2n\pi$ and $w_1(1) = 0$. 
		\item At $t=1$, change the inputs to $u_1 = \frac{\pi}{n}\sin\big(\frac{\pi}{n}t\big)$, 
		$u_2 = -\frac{\pi}{n}\cos\big(\frac{\pi}{n}t\big)$, then 
		$x_1(t) = -\cos\big(\frac{\pi}{n}t\big) + \cos\frac{\pi}{n}$, 
		$x_2(t) = -\sin\big(\frac{\pi}{n}t\big) + \sin\frac{\pi}{n}$, this realizes the curve $\gamma_a$, where 
		$a = e^{i\frac{\pi}{n}}$ in $x_1-x_2$ plane and by the above discussion at $t = 1 + 2n$, this $\gamma_a$ curve is looped around once, and $w_1(1+2n) = 2n\pi \sin(\frac{\pi}{n})$ and $w_2(1+2n) = 2n\pi - 2n\pi \cos(\frac{\pi}{n})$.
		\item Thus, the steering of this system from the origin to $(0,0,a,b)$ can be done by scaling the inputs and scaling the time taken to traverse the curves $\gamma_a$, such that $(a,b) = c_1(0,2n\pi) + c_2(2n\pi sin(\frac{\pi}{n}),- 2n\pi cos(\frac{\pi}{n}))$ 
		where $c_1$ and $c_2$ are scaling coefficients of $u_1$ and $u_2$ respectively. 
	\end{enumerate}  
	\end{example}

\begin{remark}
	For a control system given by $(\ref{fullcomplexform})$, for any continuous inputs $u_1$ and $u_2$, we cannot restrict the dynamics of 
	$(\ref{fullcomplexform})$ to the $x_1-x_2$ plane. This can be justified as follows. 
	Consider the dynamics of coordinates $w_1$ and $w_2$, we have 
	\begin{eqnarray}
	\dot{w}_1 = f_1u_1 + f_2u_2 = (f_1,f_2).(u_1,u_2) \\
	\dot{w}_2 = f_1u_2 - f_2u_1 = (-f_2,f_1).(u_1,u_2)
	\end{eqnarray}
	To restrict ourselves to $x_1-x_2$ we need to have $\dot{w}_1=\dot{w}_2=0$. Since $(f_1,f_2)$ and $(-f_2,f_1)$ are orthogonal vectors, any 
	$(u_1,u_2)$ cannot be both non-zero and be perpendicular to both of these vectors.
\end{remark}
\begin{remark}
 The complex analytic results mentioned in this section also hold for nonholonomic systems given by $(\ref{gennonhol1})$ by considering pairwise systems on the complex plane $\mathbb{C}$ 
 for all $i,j$ pairs $1\le i <j\le n$. 
\end{remark}
\section{Conclusions}
\label{sec:conclusions}

We considered generalizations of the classical nonholonomic integrator to define some specific nonholonomic systems using the notion of vector fields. 
We obtained necessary and sufficient conditions for controllability of these systems using geometric concepts such as the curl of a vector field. 
In specific, we showed that controllability is equivalent to the curl of the underlying vector field being nonzero. We also considered minimum energy 
optimal control problems on these general nonholonomic integrators and showed that the optimal trajectories are same as the trajectory of a charged 
particle in a magnetic field. We also considered a specific system with a drift term and showed that the optimal trajectories are given by a charged 
particle in an electromagnetic field. We then included a specific state dependent cost function term in the Lagrangian and showed that optimal trajectories 
are again given by the trajectory of a particle in an electromagnetic field. 
We then gave a complex analytic viewpoint to nonholonomic integrator and its generalizations and use properties such as holomorphicity, Cauchy's integral 
theorem and the residue theorem from complex analysis to characterize controllability. 

The future work involves extending these ideas for more general noholonomic systems. 
 \bibliographystyle{ieeetr}        
\bibliography{references}
\end{document}